\documentclass[draft]{amsart}
\usepackage[centertags]{amsmath}
\usepackage{amsfonts,amssymb}
\usepackage{times}

\title[Some critical issues for the ``equation-free'' approach]{Some
  critical issues for the ``equation-free'' approach\\ to multiscale
  modeling} \author{Weinan E and Eric Vanden-Eijnden}

\begin{document}

\begin{abstract}
  The ``equation-free'' approach has been proposed in recent years as
  a general framework for developing multiscale methods to
  efficiently capture the macroscale behavior of a system using only
  the microscale models.  In this paper, we take a close look at some
  of the algorithms proposed under the ``equation-free'' umbrella, the
  projective integrators and the patch dynamics.  We discuss some very
  simple examples in the context of the ``equation-free''
  approach. These examples seem to indicate that while its general
  philosophy is quite attractive and indeed similar to many other
  approaches in concurrent multiscale modeling, there are severe
  limitations to the specific implementation proposed by the
  equation-free approach.
\end{abstract}

\maketitle

\section{Introduction}
\label{sec:inro}

The purpose of this note is to examine some of the basic issues
surrounding the ``equa\-tion\--free'' approach proposed
in~\cite{eqfree}, which has been pursued in recent years as a general
tool for multiscale, multi-physics modeling.  To begin with, the
equa\-tion\--free approach is an example of concurrent coupling
technique.  In contrast to sequential coupling techniques which
require establishing the macroscale equations through precomputing,
concurrent coupling techniques compute the required macroscale
quantities ``on-the-fly'' from microscopic models \cite{AB1, AB2}.
The most well-known example of such concurrent coupling techniques is
perhaps the Car-Parrinello molecular dynamics which computes the
atomic interaction forces ``on-the-fly'' by solving the electronic
structure problem~\cite{cpmd}.  Other algorithms, such as the extended
multi-grid method \cite{Brandt} and the heterogeneous multiscale
method (HMM) \cite{hmm_cms} are all example of the concurrent coupling
approach.

At a technical level, a key idea in the ``equation-free'' approach is
to make use of scale separation in the system.  There are many
different ways of exploiting scale separation.  In~\cite{chorin} and
\cite{cpmd}, time scale separation was used to artificially slow down
the time scale of the microscopic system.  As for spatial scales,
homogenization-based methods (such as the ones that use representative
averaging volumes \cite{Bear}) and the quasicontinuum methods
\cite{QC} are all examples of algorithms that explore the separation
of spatial scales.  Most closely related to the ``equation-free''
approach is perhaps the extended multi-grid method \cite{Brandt}. In
his review article \cite{Brandt}, Achi Brandt described ideas that can
be used to extend multi-grid techniques to deal with multiscale,
multi-physics problems in order to capture the macroscale behavior of
a system using microscopic models such as molecular dynamics.  As is
common in multi-grid methods, the ideas of Brandt rely heavily on
mapping back and forth between the macro- and micro-states of the
system, through prolongation and restriction operators (which are
called respectively reconstruction and compression operators in HMM,
and lifting and restriction operators in the ``equation-free''
approach).  Brandt realized that central to the efficiency of these
algorithms is the possibility of only performing microscopic
simulations in small samples for short periods of times, as a result
of the scale separation in the system.  These ideas were later adopted
by both HMM and the ``equation-free'' approach.  In fact, HMM and the
``equation-free'' approach are both alternative approaches with the
same motivation and objective.

\vspace{.1in}
\begin{center}
  \begin{tabular}{l|l|l}
    \hline
    {} & {Macro to micro} & {micro to Macro} \\
    \hline
    Extended multi-grid &  interpolation &  restriction\\
    \hline
    HMM &  reconstruction & compression \\
    \hline
    Equation-free & lifting  & restriction \\
    \hline
  \end{tabular}
\end{center}

\vspace{.1in} While the general philosophy of the ``equation-free''
approach is very similar to the extended multi-grid method and HMM,
the ``equation-free'' approach proposes its own ways of implementing
such a philosophy, in particular, ways of dealing with scale
separation.  The basic idea is to use extrapolation in time and
interpolation in space.  More precisely, two important building blocks
of the ``equation-free'' approach are:
\begin{enumerate}
\item \underline{Projective integrators:} (An ensemble of) the
  microscale problems are solved for a short period of time using
  small time steps.  The time derivative of the macro variable is
  computed from the results of the last few steps and then used to
  advance the macro variable over a macro time step.  It is easy to
  see that such a procedure amounts to extrapolation, and indeed the
  authors state in~\cite{geke}: \textit{``The reader might think that
    these should be called `extrapolation methods,' but that name has
    already been used [...].  Hence we call the proposed methods
    projective integration methods.''}
  
\item \underline{The gap-tooth scheme:} The microscopic problem is
  solved in small domains (the teeth) separated by large gaps. The
  solution is averaged over each domain and then interpolated to give
  the prediction over the gaps.
\end{enumerate}
The combination of these two ideas gives directly the so-called
``patch dynamics''~\cite{eqfree}.

Detailed understanding of the ``equation-free'' algorithms is made
difficult by the fact that the ``equation-free'' papers are generally
quite vague.  The present note should be regarded as an attempt to pin
down some of these details.  Indeed this was initially intended as a
regular journal article.  But it soon becomes clear that there is
still substantial disagreement between our understanding of the
``equation-free'' approach and that of its developers. However, we
believe the simple examples that we discuss here do shed some light on
the ``equation-free'' approach and should be made available to a
larger audience in some form.  We are grateful to Yannis Kevrekidis
for a detailed report on the earlier version of this note. Some of his
comments have been taken into account in this revised version.  We
also welcome any discussion about the issues raised in this note, the
most important of which being: What really is the ``equation-free''
approach?  Indeed our primary purpose of presenting this note is to
prompt such a discussion.

\section{Projective integrators for stochastic ODEs}
\label{sec:projective}

Projective integrators were proposed as a way of extrapolating the
solution of an explicit ODE solver for systems with multiple time
scales using large time steps.  The basic idea is to run the
microscopic solver (using small time steps) for a number of steps, and
then estimate the time derivative and use that to extrapolate the
solution over a large time step \cite{geke}. For stiff ODEs, the
extrapolation step is applied to the whole system~\cite{geke}.  For
general multiscale problems, the extrapolation step is applied only to
the slow variables~\cite{eqfree,hummer}.

In the case of stiff ODEs, projective integrators can give rise to
useful numerical sche\-mes, as was demonstrated in~\cite{geke}.  In this
case, the idea becomes very close to the ones proposed by Eriksson
\textit{et. al} for developing explicit stiff ODE solvers~\cite{EJL}.
The objectives of the two papers are quite different: For Eriksson
\textit{et al.}, the objective is to find explicit and efficient stiff
ODE solvers.  For Gear~\textit{et al.}, the objective is to deal with
general multiscale, multi-physics problems.  However, in the general
case such as the case considered in \cite{hummer}, projective
integrators have serious limitations, as we now show.

Denote by $x$ the coarse variable of the system.  The coarse
projective integrator proposed in \cite{hummer} performs the
following steps at each macro time step (of size~$\Delta t$):
\begin{enumerate}
\item Create an ensemble of~$N$ microscopic initial conditions
  consistent with the known coarse variable~$x^n$ at time step~$n$.
\item Run the microscopic solver with these initial conditions for a
  number of steps, say~$k$, with time step~$\delta t$.  Denote the
  corresponding values of the coarse variables as~$\tilde x_j(k \delta
  t)$ where $j = 1, \ldots, N$. 
\item Perform ensemble averaging to get an approximation to the coarse
  variable. For example,
\begin{equation}
  \bar x = \frac 1N \sum_{j=1}^N \tilde x_j (k\delta t)
\end{equation}
\item Use this value to extrapolate the coarse variable to a time step
  of size~$\Delta t$.
\begin{equation}
  \label{eq:extrapolsdeb}
  x^{n+1} = x^n + {\Delta t} \frac{\bar x - x^n}{k\delta t}
\end{equation}
\end{enumerate}

Now consider the simple case when the coarse variable obeys effectively
a stochastic ODE:
\begin{equation}
  \label{eq:sdefast}
  dx(t) = b(x(t)) dt +  dW(t)
\end{equation}
Since, to  $O(k\delta t)$, we have
\begin{equation}
  \label{eq:EMstep}
  \tilde x_j(k\delta t) - x^n =  k \delta t\, b(x^n) +
  \sqrt{k\delta t} \,  \xi^n_j
\end{equation}
where $\{\xi_j^n\}, j=1, \cdots, N$ are $N$ independent Gaussian variables 
with mean 0 and variance 1, (\ref{eq:extrapolsdeb}) becomes, to leading order
\begin{equation}
  \label{eq:extrapolsde3}
  x^{n+1} = x^n + \Delta t\,  b(x^n)  + \frac{\Delta t }{\sqrt{k\delta t}}\,
  \frac1N \sum_{j=1}^N \xi_j^n.
\end{equation}
(\ref{eq:extrapolsde3}) is equivalent in law to
\begin{equation}
  \label{eq:extrapolsde4}
  x^{n+1} = x^n + \Delta t\,  b(x^n)  + 
  \frac{\Delta t }{\sqrt{N k\delta t}} \, \xi^n,
\end{equation}
where $\xi^n$ is a Gaussian variable with mean 0 and
variance 1.

It is obvious from this that the effective dynamics produced by the
coarse projective integrator depends on the numerical parameters $N$,
$k$, $\delta t, $ and $\Delta t$.  In particular, if $N k\delta t \gg
\Delta t$, then the noise term in \eqref{eq:sdefast} is lost in the
limit.  If $N k\delta t \ll \Delta t$, then the noise term overwhelms
the drift term.  In either case, one obtains a wrong prediction for
the effective dynamics of the coarse variable. 

The only way to get a scheme consistent with~(\ref{eq:sdefast}) is to
choose the numerical parameters so that they precisely satisfy $N
k\delta t = \Delta t$. The reader should be aware, however, that this
choice is not advocated in~\cite{hummer} and is, in fact, quite
orthogonal to the original equation-free philosophy since it requires
knowing beforehand that (\ref{eq:sdefast}) is an SODE and not
something else.  In addition, it is easy to see that using $N k\delta
t = \Delta t$ leads to no gain in efficiency: The total cost is
comparable to solving the microscopic problem in a brute force fashion
using $\delta t$ as the time step, since the size of the ensemble is
equal to the number of microscopic simulation time intervals
during a time duration of $\Delta t$: $ N = \Delta t/(k \delta t)$.

For the case when $N k\delta t \gg \Delta t$, one might think of using
the coarse projective integrators (or coarse molecular dynamics) as a
way of simulating the dynamics $dx/dt = b(x)$ in the context of
molecular dynamics simulations. In this case the unknown drift $b(x)$
is related to the gradient of the free energy and simulating $dx/dt =
b(x)$ is then a way to explore this free energy.  Indeed, this appears
to be how the scheme was actually used in \cite{hummer}. The problem,
however, is that using $N k\delta t \gg \Delta t$ leads again to a
scheme which is no less expensive than a brute force solution of $N$
replica of~(\ref{eq:sdefast}) using $\delta t$ as the time step.

The problem above seems to be intrinsic to projective integrators in
the context of SDEs because it is inherent to the fact that the
dynamics~(\ref{eq:sdefast}) is dominated by the noise on short time
scales and the extrapolation step in the projective integrators
amplifies these fluctuations. Averaging them out can only be done at a
cost which is at least comparable to the cost of a direct scheme.

\section{Patch dynamics}

Patch dynamics is proposed as a way of analyzing the macroscopic
dynamics of a system using microscopic models. Like the extended
multi-grid methods \cite{Brandt} and HMM \cite{hmm_cms, hmm-review},
it is formulated in such a way that scale separation can be exploited
to reduce computational cost.

The setup is as follows.  We have a macroscale grid over the
computational domain.  The grid size $\Delta x$ is chosen to resolve
the macroscale variations but not the microscale features in the
problem.  Each grid point is surrounded by a small domain (the
``tooth''), the size of which (denoted by $h$) should be large enough
to sample the local microscale variations but can be much smaller than
the macroscale grid size if the macro and microscales are very much
separated.

Given a set of macroscale values at the macroscale grid points,
$\{U^n_j\}$, at the $n$-th time step $t_n = n \Delta t$ where $\Delta
t$ is the size of the macroscale time step, patch dynamics computes
the update of these values at the next macroscale time step,
$\{U^{n+1}_j \}$, using the following procedure:

\begin{enumerate}
\item \textit{Lifting:} From $\{U^n_j\}$, reconstruct a consistent
  microscopic initial data, denoted by $\tilde u_0$.
\item \textit{Evolution:} Solve the original microscopic model with
  this initial data $\tilde u_0$ over the small domains (the
  ``teeth'') for some time $\delta t$: $\tilde u_{\delta t} =
  \mathcal{S}_{\delta t} \tilde u_0$.
\item \textit{Restriction:} Average the microscale solution $\tilde
  u_{\delta t}$ over the small domains.  The results are denoted by
  $\{ \tilde U^n_{\delta t} \}$.
\item \textit{Extrapolation:} Compute the approximate derivative and
  use it to predict $\{U^{n+1}_j \}$:
\begin{equation}
  \label{eq:extrapol}
  U^{n+1}  = U^n + \Delta t \frac{\tilde U^n_{\delta t} - U^n}{\delta t}
\end{equation}
or more generally:
\begin{equation}
U^{n+1}  = U^n + \Delta t \frac{\tilde U^n_{\delta t} - \tilde U^n_{\alpha 
\delta t}} {(1-\alpha)\delta t}
\end{equation}
where $ \alpha$ is some numerical parameter between 0 and 1.
\end{enumerate}

There are very few examples of how to implement these steps in practice.
\cite{eqfree, SKD} suggest the following:

For the lifting operator, in the small domain around the macro grid
point $x_j$, use the approximate Taylor expansion:
\begin{equation}
\label{interp-1}
\tilde u_0 (x) = \sum_{k=0}^d \frac 1{k!} D_k (x- x_j)^k
\end{equation}
Here $D_k$ is some approximations to the derivatives of the macroscale
profile at $x_j$, for example:
\begin{equation}
\label{interp-2}
D_2 = \frac{U^n_{j+1} - 2 U^n_j + U^n_{j-1}}{\Delta x^2}, \quad
D_1 = \frac{U^n_{j+1} -  U^n_{j-1}}{ 2 \Delta x}, \quad
D_0 = U^n_j - \frac 1{24} h^2 D_2
\end{equation}
Below  we will consider the case when $d=2$.

When solving the microscale problem, \cite{SKD} suggest extending the
microscale domain to include some buffer regions in the hope that this
would allow the use of {\it any} boundary conditions for the
microscopic solver: By choosing sufficiently large buffers, the effect
would be as if the microscale problem is solved in the whole space
where and when averaging is performed.  This introduces another
spatial scale $H$ which is the real size of the region on which
microscale problems are solved.  (The parameter $h$, which is (much)
smaller than $H$, is the size of the domain over which averaging is
performed).  In the following discussion, we will take $H$ to be
infinity.

Let us now examine this algorithm in more detail, using some very
simple examples.  Let us first consider the heat equation
\begin{equation}
 \partial_t u = \partial^2_{x} u 
\end{equation}
For simplicity, let us assume $x_j = 0$.
Denote $\tilde u_0 = D_0 + D_1 x + \frac 12 D_2 x^2$.
We have
\begin{equation}
\mathcal{S}_{\delta t} \tilde u_0 (x) = 
D_0 +D_1 x + D_2 (\frac 12 x^2 + \delta t)
\end{equation}
Denote by $\mathcal{A}_h$ the averaging operator over the small domain
(of size $h$), we have
\begin{equation}
\tilde U^n_{\delta t}= \mathcal{A}_h \mathcal{S}_{\delta t} \tilde u_0 (x) =
D_0 + D_2 \delta t +\frac 1{24} D_2 h^2
=U^n +D_2 \delta t
\end{equation}
Inserting this expression in~(\ref{eq:extrapol}) gives the familiar scheme:
\begin{equation}
U^{n+1} = U^n + \Delta t D_2
\end{equation}
as was shown in \cite{eqfree}.
This is both stable and consistent with the heat equation,
which is the right effective model at the large scale.

Now let us turn to the advection equation
\begin{equation}
 \partial_t u + \partial_{x} u  = 0
\end{equation}
In this case, we have
\begin{equation}
\mathcal{S}_{\delta t} \tilde u_0 (x) = 
D_0 +D_1 (x -\delta t)  + \frac 12 D_2 (x - \delta t)^2 
\end{equation}
Hence, 
\begin{equation}
  \tilde U^n_{\delta t}= \mathcal{A}_h \mathcal{S}_{\delta t}  \tilde u_0 (x)=
  D_0 - D_1 \delta t + 
  \frac 12 D_2 \delta t^2 +\frac 1{24} D_2 h^2 
  = U^n -D_1 \delta t + \frac 12 D_2 \delta t^2
\end{equation}
and (\ref{eq:extrapol}) becomes
\begin{equation}
U^{n+1} = U^n  +\Delta t (-D_1+ \frac 12 \delta t D_2)
\end{equation}
Since $\delta t \ll \Delta t$, the last term is much smaller than the other
terms, and we are left essentially with a scheme which is unstable under
the standard CFL condition that $\Delta t \sim \Delta x$:
\begin{equation}
U^{n+1} = U^n  -\Delta t D_1
\end{equation}
due to the central character of $D_1$.

Aside from the stability issue, there can also be problems with
consistency. Consider the following example:
\begin{equation}
\partial_t u = - \partial_x^4 u
\end{equation}
The macroscale model is obviously the same model.  However, it is easy
to see that if we follow the patch dynamics procedure with $d=2$, we
would be solving $\partial_t U = 0$, which is obviously inconsistent
with the correct macroscale model.

For these simple examples, the difficulties discussed above can be
fixed by using different reinitialization procedure for the
micro-solvers.  For the example of the convection equation, one should
use one-sided interpolation schemes in the spirit of upwind
schemes. For the last example, one should use piecewise $4th$ order
polynomial interpolation.  {\it But in general, finding such a
  reinitialization procedure seems to be quite a daunting task, since
  it depends on the nature of the unknown effective macroscale model}.
Imagine that the microscopic solver is molecular dynamics.  The
reinitialization procedure has to take into account not only
consistency with the local macrostates of the system (which is the
only requirement for the extended multi-grid method and HMM), but also
the effective macroscale scale model (which is unknown) such as:
\begin{enumerate}
\item The order of the macroscale equation.
\item The direction of the wind, if the effective macroscale 
equation turns out to be a first order PDE.
\item Other unforeseeable factors.
\end{enumerate}
Indeed it is not at all clear how patch dynamics would work if
molecular dynamics models are used to model macroscopic gas dynamics.

To overcome these problems, the ``equation-free'' developers propose
to design a number of numerical tests to find out the nature of the
effective macroscale equations.  One such a procedure, the
``baby-bathwater scheme'' will be discussed in the next section.
However, in addition to  the technical issues, it is not clear what set of
characteristics that we are supposed to test on.

\section{The ``baby-bathwater scheme''}

As the last example shows, it is useful at times to know the order of
the highest order derivatives that appear in the effective macroscale
equation, even if we do not know all the details of the macroscale
model.  An algorithm was proposed in \cite{baby} for this purpose.
The ``baby-bathwater scheme'', as it was called, promises to find the
highest order derivative in the effective macroscale model, by
performing simulations using the microscopic model: Assume that the
effective macroscale model is of the form
\begin{equation}
\partial_t U = F(U, \partial_x U, \cdots, \partial_x^m U)
\end{equation}
The objective is to find $m$.

This problem can be formulated abstractly as follows.  Assume we have
a function $F=F(x_1, x_2, \cdots)$ and we know that it only depends on
finitely many variables: $F= F(x_1, x_2 \cdots, x_m)$.  Assume that we
can evaluate $F$ at any given point, can we find the value of $m$
efficiently?

The basis idea used in \cite{baby} is that if $F$ depends truly
on the variable $x_j$, then the variance of $F$ as $x_j$ changes should
not vanish. The practical difficulty is how to use this idea
efficiently.

Without getting into the details of the algorithm presented in
\cite{baby}, it seems quite clear that there is no fool-proof
inductive procedure for finding $m$. Take an extreme case, say, $F=
F(x_1, x_{100})$.  Without having the prior knowledge that $F$ may depend
on $x_{100}$, an inductive procedure would likely conclude that $F$ is
only a function of $x_1$.

This example is of course quite extreme, and most practical situations
are not like this. Nevertheless, it does raise some questions about
the robustness of the algorithm presented in \cite{baby}.  There is a
much more serious concern, and this is associated with the well-known
phenomenon that the order of the effective macroscale model depends on
the scale we look at.  A physically intuitive example is convection
and diffusion of tracer particles in the Rayleigh-Bernard cells: At the
scale of the cells, the tracer particles are convected and diffusion
can be neglected.  Hence the effective model is a first order
equation.  At scales much larger than the size of the cells, diffusion
dominates.  Hence the effective model is a second order equation.
This means that the output of the ``baby-bathwater scheme'' should
depend on the numerical parameter $\Delta$ in the scheme.

This behavior can be demonstrated rigorously using the well-known
results of Kesten and Papanicolaou \cite{K-P}.  Consider the dynamics
if inertial particles in a stationary random force field:
\begin{equation}
  \frac{d^2 x}{d t^2} = F(x)
\end{equation}
In phase space, we can write this as
\begin{equation}
  \begin{cases}
    \displaystyle \frac{d x}{d t} = v,\\[6pt]
    \displaystyle \frac{d v}{d t} = F(x)
  \end{cases}
\end{equation}
The density of the particles (in phase space)
obeys the Liouville equation:
\begin{equation}
\partial_t \rho + v \partial_x \rho + F(x) \partial_v \rho = 0
\end{equation}
However, if we consider the rescaled fields:
\begin{equation}
  \begin{cases}
    \displaystyle \frac{d x^\delta}{d t} = \frac1{\delta^2} v^\delta,\\[6pt]
    \displaystyle \frac{d v^\delta}{d t} = \frac1{\delta} F(x^\delta)
  \end{cases}
\end{equation}
it was shown in \cite{K-P} under some conditions on the random field
$F$ that as $\delta \rightarrow 0$, the process $v^\delta (\cdot)$
converges to a diffusion process.  In other words, if we consider the
density of the particles in $v$ space, then in this scaling we have
\begin{equation}
  \partial_t \rho +  \partial_v 
  (b(v) \rho)   =  \frac 12 \partial_v^2 (a(v) \rho)
\end{equation}
for some functions $b(\cdot)$ and $a(\cdot)$, which is a second order
equation.

\section{Conclusions}

The idea of interrogating legacy codes as a control system is very
attractive and to some extend, has already been commonly used in some
disciplines.  For example, chemists use packages such as CHARMM and
AMBER as legacy codes to perform optimization tasks, e.g.  to find
free energy surfaces and minimum free energy paths.  Optimization
techniques such as the Nelder-Mead algorithm that use only function
values (not the derivatives) were designed with this kind of problems
in mind.  One purpose of the work of Keller \textit{et al.} is to
extend bifurcation analysis tools to systems that are defined by
legacy codes \cite{Shroff}.  The ``equation-free'' approach attempts
to extend such practices to another direction, namely the modeling of
macroscale spatial/temporal dynamics of systems defined by microscopic
models, in the form of legacy codes.  While this seems very
attractive, the set of tools proposed under this umbrella are quite
far from being sufficient for reaching this objective.  We have
discussed some of the technical difficulties in this note.  This
discussion is certainly not exhaustive.  It is only meant to be
illustrative.

\section*{Acknowledgments}
We are very grateful to the many people, including Yannis Kevrekidis,
who have read the draft of this paper at various stages and gave us
their comments.  The work of W. E was supported in part by ONR grant
N00014-01-1-0674 and by DOE grant DOE DE-FG02-03ER25587.  The work of
E.  V.-E. was supported in part by NSF grants DMS02-09959 and
DMS02-39625 and by ONR grant N00014-04-1-0565.

\bibliographystyle{elsart-num}

\end{document}